\documentclass[12pt,a4paper,twoside]{article}

\usepackage[T1]{fontenc}
\usepackage{amsmath}
\usepackage{amssymb}

\usepackage{graphicx} 
\usepackage{psfrag}

\def\newtext1 {}

\def\0{{\bf {0}}}

\makeatletter
\renewcommand{\subsection}{\@startsection{subsection}{2}{0mm}{\baselineskip}%
{-\fontdimen2\font plus -\fontdimen3\font minus \fontdimen4\font}%
{\normalfont\normalsize\bfseries}}
\renewcommand{\@seccntformat}[1]{\csname the#1\endcsname .\hspace{0.5em}}
\makeatother

\newtheorem {lemma}{Lemma}

\newtheorem {theorem}[lemma]{Theorem}
\newtheorem {proposition}[lemma]{Theorem}

\newtheorem {corollary}[lemma]{Corollary}

\def\vol{\hbox{vol}}
\def\M{{M}}

\def\tilde{\widetilde}

\def \mbeq {\begin {eqnarray}}
\def \meeq {\end {eqnarray}}

\def \HL {\Lambda_{2T}}

\def \bfo {\begin {displaymath} }
\def \efo {\end {displaymath} }

\def \beq {\begin {eqnarray}}
\def \eeq {\end {eqnarray}}
\def \ba {\begin {eqnarray*}}
\def \ea {\end  {eqnarray*}}

\def \Z {{\Bbb {Z}}}

\def \R {{\Bbb {R}}}
\def \expec {{\Bbb {E}}}
\def \prob {{\Bbb {P}}}

\def \bx {{{x}}}

\def \H2s {H^{s+1}_0(\partial \M\times [0,T/2])}

\def \supp {\hbox{supp }}

\def \dist {d}

\def \det {\hbox{det}}
\def\bra{\langle}
\def\cet{\rangle}

\def \e {\varepsilon}

\def \a {\alpha}

\def \pa0 {\partial _0}

\def \p {\partial}

\def \tilde{\widetilde}
\def \proofbox {$\Box$\medskip}
\def \e {\epsilon}

\begin{document}

\title{
Time reversal methods in unknown medium and inverse problems
}
 \author{Kenrick Bingham\\
 Yaroslav Kurylev\\
 Matti Lassas\\ Samuli Siltanen
}
\date{January 4, 2007
}

\maketitle

\begin{abstract}
A novel method to solve inverse problems for the wave equation is introduced. 
The method is a combination of the boundary control method and an iterative
time reversal scheme, leading to adaptive imaging of coefficient functions
of the wave equation using focusing waves in unknown medium. 
The approach is
computationally effective since the iteration lets the medium do 
most of the processing of the data.

The iterative time reversal scheme also gives an algorithm for
constructing boundary controls for which the corresponding final values 
are as close as possible to the final values of a given wave 
in a part of the domain, 
and as close as possible to zero elsewhere. The algorithm does not
assume that the coefficients of the wave equation are known.
\end{abstract}



\section{Introduction}

We present a novel inversion method for the wave equation. Suppose that we can send waves from the boundary into an unknown body with
spatially varying wave speed $c(x)$. Using a combination of the boundary control (BC) method and an iterative time reversal scheme, we show
how to focus waves near a point $x_0$ inside the medium and simultaneously recover $c(x_0)$ if $c$ is isotropic. In the anisotropic case we
can reconstruct $c(x)$ up to a change of coordinates.

In the classical BC method material parameters are reconstructed using hyperbolic techniques, see \cite{AKKLT,B1,B2,BK,KK,KKL,KKLima}.
Straightforward numerical implementation of the BC method is computationally demanding. Inspired by Isaacson's iterative
measurement scheme \cite{I} for electrical impedance tomography (EIT), we combine the BC method with an adaptive time
reversal iteration that lets the medium do most processing of the data.

Traditional time reversal methods record waves, invert them in time, and send them back into the medium. If the recorded signal originated
from point sources or was reflected from small scatterers in the medium, the time reversed waves focus at the source or scatterer points.
This is useful for time reversal mirrors in communication technologies and medical therapies. For early references on time reversal (using
ultrasound in air), see the seminal works of Fink \cite{FinkMain,FinkB,FinkD}. For a microlocal discussion of time reversal, see Bardos
\cite{Bardos,Bardos2}, and for another mathematical treatment see Klibanov \cite{Kliba}.

Time reversal in known background medium with random fluctuations has also been extensively studied \cite{Bal1,Bal2,Bal3,Bal4,Papa1,Papa3}
and applied to medical imaging, non-destructive testing and underwater acoustics. These methods are outside the scope of this paper.

Iterative time reversal has been used to find best measurements for inverse problems. By ``best'' we here mean ``optimal for detecting the
presence of an object''. This {\em distinguishability problem} has been studied for fixed-frequency problems in EIT \cite{I} and acoustic
scattering \cite{MNW}.  The connection between optimal measurements and iterative time-reversal experiments was pointed out in
\cite{MNW,PF,PTF}. For the wave equation, the best measurement problem has been studied in \cite{CIL}, where it has been shown that the
optimal incident field for probing a half space can be found by an iterative process involving time reversal mirrors.

Let us describe our hybrid method in a simple case. Take a compact set $M\subset \R^3$ with smooth boundary $\partial M$, and let $c(x)$ be
a scalar-valued wave speed in $M$. Consider the wave equation
\begin{align}\label{eq: Wave isotropic}
  u_{tt}-c(x)^2\Delta u&=0 \qquad \hbox{ in } M\times \R_+,
  \\
  u|_{t=0}&=0,\nonumber\\
  u_t|_{t=0}&=0,  \nonumber \\
  -c(x)^{-2}\p_n u&= f(x,t) \qquad \hbox{ in }\p M\times \R_+,\nonumber
\end{align}
where $\p_n$ denotes the Euclidean normal derivative. We denote by $u^f=u^f(x,t)$ the solution of (\ref{eq: Wave isotropic}) corresponding
to a given boundary source term $f$.

Take $T>0$ larger than the diameter of $M$ in travel time metric, so that any point inside $M$ can be reached by waves sent from the
boundary in time less than $T$. We consider imposing a boundary source $f$ on the set $\p M\times (0,2T)$ and measuring the boundary value
of the resulting wave $u^f$ during the time $0 < t < 2T$. All such boundary measurements constitute the Neumann-to-Dirichlet map
$$
  \Lambda_{2T}:f\mapsto u^f|_{\p M\times (0,2T)},
$$
that models a variety of physical measurements \cite{KKLM}. Instead of the full map we only need here the values $\Lambda_{2T}f$ for
a specific collection of sources $f$ given by an iterative process.

We need three special operators on the function space $L^2(\p M\times [0,2T])$. First, the {\em time reversal operator}
$$
 Rf(x,t)=f(x,2T-t),
$$
second, the {\em time filter}
$$
 Jf(x,t)= \frac 12 \int_0^{\min(2T-t,t)}f(x,s)ds,
$$
and finally the {\em restriction operator}
$$
  Pf(x,t)=\chi_B(x,t)u(x,t),
$$
where $\chi_B$ is the characteristic function of a set $B\subset \p M\times [0,2T]$.

We define the \emph{processed time reversal iteration}
\begin{align} \label{eq: iteration 1 motivate}
F&:= \frac 1\omega P(R \Lambda_{2T} R J-J\Lambda_{2T})f,\nonumber\\
a_n&:=\Lambda_{2T}(h_n),\qquad
b_n:=\Lambda_{2T}(R J h_n),\\
\nonumber h_{n+1}&:=(1-\frac \alpha \omega)h_n- \frac 1\omega(PR b_n-PJ a_n) + F,
\end{align}
where $f\in L^2(\p M\times [0,2T])$ and $\alpha,\omega>0$ are parameters. Starting with $h_0=0$ the iteration (\ref{eq: iteration 1
motivate}) converges to a limit $h(\alpha)=\lim_{n\rightarrow\infty}h_n$.
In the case $B=\p M\times [T-s,T]$ we have
\begin{equation}\label{eq: limits}
  \lim_{\alpha\to 0} u^{h(\alpha)}(x,T)= \big(1-\chi_N(x) \big) u^f(x,T)\quad \hbox{in }L^2(M),
\end{equation}
where $N\subset M$ is the set of points $x\in M$ with travel time distance to $\p M$ at most $s$. In other words, $N$ is the maximal region
where the waves sent from the boundary can propagate in time $s$.

Let $y\in \p M$ and denote by $\gamma_{y,\nu}$ the geodesic (in travel time metric) starting from $y$ in the normal direction and
parametrised by arclength. Take a source $f$ concentrated near $(y,T-s)\in \p M\times \R_+$ and vanishing for $t\leq T-s$. Then $u^f(x,T)$
is supported in the union of the set $N$ and a neighbourhood of the point $x_0=\gamma_{y,\nu}(s)$, and the limit (\ref{eq: limits})
vanishes outside a neighbourhood of $x_0$. This way our iteration produces a wave localized near $x_0$ at time $T$.

Below we introduce an alternative construction of focusing waves for general sources $f$ by applying the processed time reversal iteration
with two different sets $B$ and using the difference of the resulting boundary values.

We show also that the limits (\ref{eq: limits}) can be used to determine the travel time distance between an arbitrary point
$x=\gamma_{y,\nu}(s)\in M$ and any boundary point $z\in \p M$. It follows that the processed time reversal iteration gives an algorithm for
determining the wave speed in the medium. Furthermore, given a source $f$, our iteration allows us to find a wave $u^h(x,t)$ having at time
$t=T_1$ approximately the value $\chi_N(x)u^f(x,T_2)$ for any $T_2$ without knowing the material parameters of the medium. Analogous
methods without the multiplier $\chi_N(x)$ have been developed before in \cite{Jon-DeHoop,KL_edin}.

The paper is organised as follows. In Section \ref{sec: results} we formulate our main results. In Section \ref{sec:localization proofs} we
prove the convergence of the processed time reversal iteration. Section \ref{sec:focus} is devoted to the analysis of the focusing
properties of the waves, and in Section \ref{sec:metric} we show how to reconstruct the metric using the boundary distance function. In the
last section we discuss our method in the case of noisy measurements.

\section{Definitions and main results} \label{sec: results}

Let us consider the closure
$M\subset \R^m$, $m\geq 1$,
of an open smooth set, or a  (non-compact or compact) complete
Riemannian manifold $(M,g)$ of dimension $m$ with a non-empty boundary.
For simplicity, we assume that the boundary $\p M$ is compact.
 Let $u$
solve the  wave equation
\begin{align}\label{eq: Wave}
  u_{tt}(x,t)+Au(x,t)&=0\quad \hbox{ in }\quad M\times \R_+,\\
u|_{t=0}&=0,\quad u_t|_{t=0}=0,  \nonumber \\
B_{\nu,\eta}u|_{\p M\times \R_+}&=f. \nonumber
\end{align}
Here, $f\in  L^2(\p M\times \R_+)$ is a real valued function,
$A$ is a formally self-adjoint elliptic partial differential
operator of the form (in local coordinates in the case when $M$ is a manifold)
\beq\label{operator A}
Av=-\sum_{j,k=1}^m \mu(x)^{-1}|g(x)|^{-\frac 12}\frac {\p}{\p x^j}\left(
\mu(x)|g(x)|^{\frac 12}g^{jk}(x)\frac {\p v}{\p x^k}(x)\right)+q(x)v(x)
\eeq
where $[g^{jk}]$ is a smooth real positive definite matrix,
$c_0I\leq [g^{jk}(x)]\leq c_1I$, $c_0,c_1>0$, $|g|=\det([g_{jk}])$,
where $[g_{jk}]$ is the inverse matrix of  $[g^{jk}]$,
 and $\mu(x)\geq c_2>0$ and $q(x)$
are smooth real valued functions. Also,
\ba
B_{\nu,\eta}v=-\p_\nu v+\eta v
\ea
where $\eta:\p M\to \R$ is a smooth bounded function and
\ba
\p_\nu v=\sum_{j,k=1}^m \mu(x)g^{jk}(x)\nu_k\frac {\p}{\p x^j}v(x),
\ea
where $\nu(x)=(\nu_1,\nu_2,\dots,\nu_m)$ is the interior co-normal vector
field of $\p M$ normalised so that
$\sum_{j,k=1}^m g^{jk}\nu_j\nu_k=1$.
A particular example is the operator
\beq\label{eq: example A}
A_0=-c^2(x)\Delta+q(x)
\eeq
{for which $\p_\nu v= c(x)^{-m+1} \p_n v$, where
$\p_n v$  is the Euclidean normal derivative of $v$. }

We denote the solutions of (\ref{eq: Wave}) by
\ba
u^f(x,t)=u(x,t).
\ea
For the initial boundary value problem  (\ref{eq: Wave})
we define the response operator (or non-stationary Robin-to-Dirichlet map)
$\Lambda$ by setting
\begin{equation}
\label{1.5}
\Lambda f = u^f|_{\p M\times \R_+}.
\end{equation}
We also consider the finite time response operator $\Lambda_T$ corresponding
to time $T>0$,
\begin{equation}
  \label{1.5b}
\Lambda_T f = u^f|_{\p M\times (0,T)}.
\end{equation}
By \cite{tataruREG}, the map $\Lambda_T:
L^2(\p M\times (0,T))\to H^{1/3}(\p M\times (0,T))$ is bounded,
where $H^{s}(\p M\times (0,T))$ denotes the Sobolev space on
$\p M\times (0,T)$. Below we consider $\Lambda_T$ as a bounded
operator that maps $L^2(\p M\times (0,T))$ to itself.

We also need some other notation.
The matrix $g_{jk}(x)$,  that is, the inverse of the matrix $g^{jk}(x)$,
is a Riemannian
metric in $M$ that is called the travel time metric. The
reason for this is the waves propagate with speed one
with respect to  the metric $ds^2=\sum_{jk} g_{jk}(x)dx^jdx^k$.
We denote by $d(x,y)$ the distance function
corresponding to this metric.
For the wave equation we define the space
$L^2(M,dV_\mu)$ with inner product
\ba
\bra u,v\cet_{L^2(M,dV_\mu)}=\int_M u(x)v(x)\,{dV_\mu(x)},
\ea
where $dV_\mu=\mu(x)|g(x)|^{1/2}dx^1dx^2\dots dx^m$.

For $t>0$ and $\Gamma\subset \p M$, let
\beq
\label{26.2}
M(\Gamma,t)=\{x\in M\ :\ d(x,\Gamma)\leq t\},
\eeq
be the domain of influence of $\Gamma$ at time $t$.

We denote $ L^2(B)=\{f\in
 L^2(\p M\times \R_+):\ \supp(f)\subset B\}$, $B\subset \p M\times \R_+$,
identifying functions and their zero continuations.
When  $\Gamma\subset \p M$ is open set and
  $f\in L^2(\Gamma\times \R_+)$, it is well known (see
e.g.\ \cite{H3}) that the wave $u^f(t)=u^{f}(\cdotp,t)$ is supported in the
domain $M(\Gamma,t)$,
\ba
u^{f}(t)\in L^2(M(\Gamma,t))=\{v\in L^2(M):\ \supp(v)\subset
M(\Gamma,t)\}.
\ea

\subsection{Results on convergence of processed time reversal iteration.}
Our objective is to find a boundary source $f$ such that the wave
$u^f(x,T)$ is localized in a neighbourhood of a single
point $x_0$. Note that in this paper we do not focus
the pair $(u(x,T),u_t(x,T))$, but only the value of the wave,
$u(x,T)$.

Let  $0<T_0<T$, $\Gamma\subset \p M$ be open
 set
and  $f\in L^2(\Gamma\times \R_+)$. Our first aim
is to construct boundary values $h(\alpha)\in L^2(\p M\times \R_+)$,
 $\alpha>0$ such that
\beq\label{objective eq: 1}
\lim_{\alpha\to 0} u^{h(\alpha)}(T)=\chi_{M(\p M,T_0)} u^f(T)
\eeq
in $L^2(M)$.
Here $\chi_X(x)$ is the characteristic function of the set $X$.
Then
\ba
\lim_{\alpha\to 0}u^{f-h(\alpha)}(T)=
(1-\chi_{M(\p M,T_0)} )u^f(T)
\ea
 is supported in
$M(\Gamma,T)\setminus M(\p M,T_0)$ (see Figure 1).
When $T_0\to T$ and $\Gamma_j\to \{z\}$ then under suitable
assumptions discussed later the set
$M(\Gamma_j,T)\setminus M(\p M,T_0)$ goes to a single point.
Thus  for $f\in L^2(\Gamma_j\times \R_+)$, the waves
$(1-\chi_{M(\p M,T_0)} )u^f(T)$ localize to a small neighbourhood of
a single point. For later use we will consider (\ref{objective eq: 1})
below in a more general setting.

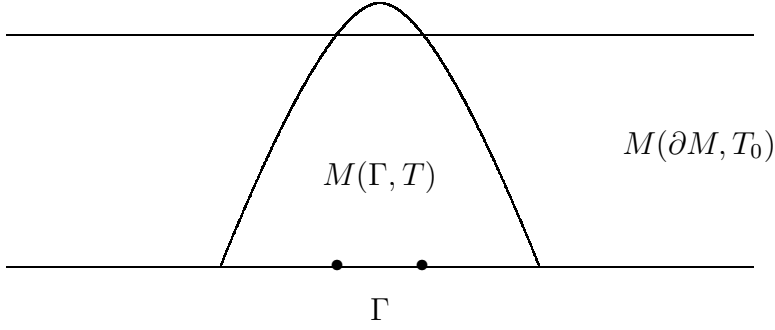
\begin {figure}
\caption{The set $M(\Gamma,T)\setminus M(\p M,T_0)$. Our goal
is to find waves sent from the boundary that focus into this area.}
\label {Fig:4.4.1a}
\setlength {\unitlength}{1.4mm}
\begin {picture}(100,45)(0,0)
\thinlines
\qbezier (30,10)(40,35)(45,35)
\qbezier (45,35)(50,35)(60,10)

\qbezier (80,10)(80,10)(45,10)
\qbezier (45,10)(10,10)(10,10)

\qbezier (80,32)(80,32)(45,32)
\qbezier (45,32)(10,32)(10,32)

\put (49,10.2){\circle*{1}}
\put (41,10.2){\circle*{1}}
\put (45,5){\makebox(0,0)[b]{$\Gamma$}}
\put (75,20){\makebox(0,0)[b]{${M(\p M,T_0)}$}}
\put (45,17){\makebox(0,0)[b]{${M(\Gamma,T)}$}}

\end {picture}
\end {figure}

Next we explain a procedure to find boundary values $h(\alpha)$.
Denote
\beq \label{eq:notations}
& &R_{2T}f(x,t)=f(x,2T-t),\\ \nonumber
& &J_{2T}h(x,t)=  \int_{[0,2T]} J(s,t)h(x,s)ds
\eeq
where $J(s,t)=\frac 12\chi_{L}(s,t)$,
\beq\label{Eq: L}
L=\{ (s,t)\in \R_+\times\R_+:\ t+s\leq 2T,\ \ s > t \}.
\eeq
We call $R=R_{2T}$ the \emph{time reversal map} and $J=J_{2T}$ the
\emph{time filter}. We denote by $P=P_B:L^2(\p M\times [0,2T])
\to L^2(\p M\times [0,2T])$
the multiplication operator
\ba
P_Bf(x,t)= \chi_B(x,t)\,f(x,t),
\ea
where $B=\bigcup_{j=1}^J (\Gamma_j\times [T-T_j,T]),$
$\Gamma_j\subset \p M$ are open sets and $0\leq T_j\leq T$.
Next, we consider $\Lambda_{2T}$, $R_{2T},$ and $J_{2T}$ as operators
from $L^2(\p M\times [0,2T])$ to itself.

Let $\alpha\in (0,1)$ and $\omega>0$
 be a sufficiently large constant.
We define $a_n,b_n\in L^2(\p M\times [0,2T])$ and
$h_n=h_n(\alpha)\in L^2(\p M\times [0,2T])$
by the iteration
\begin{align}
\label{eq: iteration 1}
a_n&:=\Lambda_{2T}(h_n),\quad
b_n:=\Lambda_{2T}(R J h_n),\\
\nonumber h_{n+1}&:=(1-\frac \alpha \omega)h_n- \frac 1\omega(PR b_n-PJ a_n) +F,
\end{align}
with $a_0=0$, $b_0=0$, $h_0=0$, and \ba F=\frac 1\omega P(R \Lambda_{2T} R J-J\Lambda_{2T})f. \ea We say that this iteration is the
\emph{processed time reversal iteration on time-interval $[0,2T]$ with projector $P$} and starting point $f\in L^2(\p M\times [0,2T])$.

Here $a_n$ corresponds to the ``iterated measurement'',
$b_n$ to the ``time filtered and time-reversed measurement''
and $h_{n+1}$ to post-processing of $h_n$, $a_n$ and $b_n$
using time reversal and time filtering.

\begin{theorem}\label{Main F} Let $T>0$.
Assume we are given $\p M$ and the response
operator $\Lambda_{2T}$. Let
 $\Gamma_j\subset \p M$, $j=1,\dots,J$
 be  non-empty open sets, $0\leq T_j\leq T$,
and  $B=\bigcup_{j=1}^J (\Gamma_j\times [T-T_j,T]).$
Let  $f\in L^2(\p M\times \R_+)$ and let
for $\omega$ large enough and $\alpha\in (0,1)$
functions $h_n=h_n(\alpha)$ be defined by
the processed time reversal iteration (\ref{eq: iteration 1})
with projector $P_B$ and starting point $f_0=f|_{\p M\times (0,2T)}$.
These functions converge in $L^2(\p M\times \R_+)$,
\ba
h(\alpha)=\lim_{n\to \infty }h_n(\alpha)
\ea
and the limits satisfy
\ba
\lim_{\alpha\to 0} u^{h(\alpha)}(x,T)=\chi_{N}(x) u^f(x,T)
\ea
in $L^2(M)$, where $N=\bigcup_{j=1}^J M(\Gamma_j,T_j)\subset M$.
\end{theorem}

This theorem is proven later in Section \ref{sec:localization proofs}.
Before that we consider its consequences.

\subsection{Results on focusing of the waves}


Let us consider a geodesic $\gamma_{x,\xi}$ in $(M,g)$
parametrised along arc length where
$\gamma_{x,\xi}(0)=x$,
$\dot\gamma_{x,\xi}(0)=\xi$ with $\|\xi\|_g=1$.
Let $\nu=\nu(z)$, $z\in \p M$ be the unit interior normal vector of $\p M$ in $(M,g)$.
There is a critical value $\tau(z)\in (0,\infty]$, such that for
$t<\tau(z)$ the geodesic $\gamma_{z,\nu}([0,t])$ is the unique
shortest geodesic, and for
$t>\tau(z)$ it is no longer a shortest geodesic from its endpoint
$\gamma_{z,\nu}(t)$ to $\p M$.


We say that $\Gamma_j \to \{\widehat{z}\}$ if $\Gamma_{j+1} \subset \Gamma_{j}$
and $\bigcap_{j=1}^\infty \Gamma_j=\{\widehat{z}\}$.

Theorem \ref{Main F} yields the following results telling that
we can produce focusing waves,
that is, wave that focus to a single point.

\begin{corollary}\label{Cor 1}
Let
${\widehat z}\in \p M$ and $0 < T_0 < \widehat T < T$.
Let  ${\widehat x}=\gamma_{{\widehat z},\nu}(\widehat T)$
and  $\Gamma_j\subset \p M$, $j\in \Z_+$ be open neighbourhoods
of ${\widehat z}\in \p M$ such that $\Gamma_j\to \{{\widehat z}\}$ when $j\to\infty$.

Let $f\in C^\infty_0(\p M \times \R_+)$. Let
$h_{n}(\alpha; T_0, j)$ be the functions obtained from
the processed time reversal iteration (\ref{eq: iteration 1})
with projector $P_B$,  $B=(\Gamma_j\times [T-\widehat T,T])\cup
(\p M\times [T-T_0,T])$
 and starting point $f_0=f|_{\p M \times (0,2T)}$. Similarly, let
$h_{n}'(\alpha; T_0, j)$ be the functions obtained from
the processed time reversal iteration (\ref{eq: iteration 1})
with
projector $P_{B'}$,
 $B'=\p M\times [T-T_0,T]$ and starting point $f_0$. Let
$$
\tilde h_n(\alpha; T_0, j)=f_0-
h_{n}(\alpha; T_0, j)+h_{n}'(\alpha; T_0, j).
$$
If $\widehat T< \tau(\widehat z)$ then
\beq\label{a limit}\\
\nonumber
\lim_{T_0\to \widehat T}\lim_{j\to \infty}\lim_{\alpha\to 0}\lim_{n\to \infty}
\frac 1{(\widehat T-T_0)^{(m+1)/2}}\,u^{\tilde h_{n}(\alpha;T_0,j)}(T)=
C_0(\widehat x)u^f(\widehat x, T)\delta_{\widehat x}(x)
\eeq
in ${\mathcal D}'(M)$ where $C_0(\widehat x)> 0$
does not depend on $f$.
If $\widehat T> \tau(\widehat z)$ the limit (\ref{a limit}) is zero.

\end{corollary}

\begin{figure}
  \caption{Subsets of $M$ used in the proof of Corollary~\ref{Cor 1}.}
  \label{fig:focusingslice}
  \centering
  \hspace{-.2cm}\raisebox{-1cm}
  \mbox{
    \psfrag{1}{$\p M$}
    \psfrag{2}[c]{$\widehat z$}
    \psfrag{3}[c]{$\Gamma_j$}
    \psfrag{4}[t]{$\widehat x$}
    \psfrag{6}{$T_0$}
    \psfrag{7}{$\widehat T$}
    \includegraphics[width=9cm]{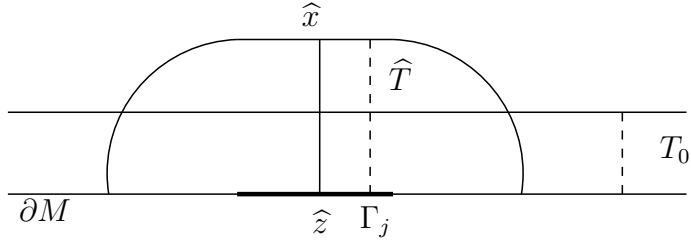}}
  \hspace{.2cm}
\end{figure}

In Figure~\ref{fig:focusingslice}, the set
$M(\Gamma_j,\widehat T)\setminus M(\p M,T_0)$
tends to the point $\hat x$ when $j\to \infty$ and $T_0\to \hat T$.
This turns out to be essential in the proof of Corollary \ref{Cor 1}.

Above, $\delta_{\widehat x}$ is Dirac's delta distribution on $(M,g)$ such that
\ba
\int_M \delta_{\widehat x}(x) \phi(x)\,dV_\mu= \phi(\widehat x),\quad\phi\in C^\infty_0(M).
\ea

\subsection{Inner products.}
Later, we show that that using
 $f,h\in L^2(\p M\times [0,2T])$ and the boundary measurements
$\Lambda_{2T}$  we may compute via an explicit integral
the inner product
\beq\label {4.60b}
\int_M  u^f(x,T)u^h(x,T)\,dV_\mu(x)=\int_{\p M\times [0,2T]} (Kf)(x,t)h(x,t)\,
dS_g(x)dt,
\eeq
where
$dS_g$ is the Riemannian surface volume of $\p M$
and
\beq\label{K-definition}
K=K_{2T}:=R_{2T}\Lambda_{2T} R_{2T} J_{2T}-J_{2T}\Lambda_{2T}
\eeq
so that terms (\ref {4.60b})
can be found using measurement operator $\Lambda_{2T}$ and
simple basic operations like time reversal $R_{2T}$ and the
time filter operator $J_{2T}$.

In the above formula the Riemannian surface volume of $\p M$
can be computed using the  intrinsic metric of $\p M$, which is  determined
by the map $\Lambda_{2T}$. Indeed,
 its follows
from Tataru's unique continuation principle  (see \cite{Ta1,Ta3},
see also Theorem \ref{th:3.4} below) that
the Schwartz kernel of $\Lambda_{2T}$ is supported
in the set $E=\{(x,t,x',t')\in (\p M\times [0,2T])^2:\ t-t'\geq
d(x,x')\}$ and that the boundary $\p E$ is in the support.
The set $\p E$ determines the distances of points $z,z'\in \p M$
with respect to the metric of $(M,g)$, and thus also
the distance
with respect to the intrinsic metric of the boundary $(\p M,g_{\p M})$.

\subsection{Reconstruction of material parameter functions.}

It is well known that the response map $\Lambda$ can
not generally determine coefficients $\mu$ and $g_{jk}$
because of two transformations discussed below.

 First, we one can introduce
a coordinate transformation, that is a diffeomorphism $F:M\to M$ such that the boundary value
$F|_{\p M}$ is the identity operator. Then
the push forward metric $\tilde g=F_*g$, that is,
\ba
\tilde g_{jk}(y)=\sum_{p,q=1}^m\frac {\p x^p}{\p y^j}\frac {\p x^q}{\p y^k}\, g_{pq}(x),
\quad y=F(x)
\ea
and the functions
$\tilde \mu=\mu\circ F^{-1}$, $\tilde q=q\circ F^{-1}$, and
$\tilde \eta=\eta$
determine the operator $\tilde A$ of the form (\ref{operator A})
such that response operators
for $A$ and $\tilde A$ are the same (for this, see \cite{KK, KKL,KKLima}).

Second, one can do the
gauge transformation $u(x,t)\to \kappa(x)u(x,t)$
where
$\kappa\in C^{\infty}(M)$ is strictly positive function.
In this transformation
the operator $\tilde A$ is transforms to the operator
$\tilde A_\kappa$ defined be
\ba
\tilde A_\kappa w:=\kappa \tilde A(\frac 1\kappa w)
\ea
and the operator in boundary condition is transformed to $B_{\nu,\widehat \eta}$,
with
\ba
\widehat  \eta=\eta-\kappa^{-1}\p_\nu \kappa.
\ea
Then the  response operator $\tilde \Lambda_\kappa$ of
$\tilde A_\kappa$ coincides with the response
operator $\Lambda$ of $A$.

By above transformations, making a gauge transformation with $\kappa=\mu^{-1}$,
we come to operator $\tilde A_\kappa$ having
the same response operator as $A$ and that can be
represented in the form (\ref{operator A}) with $\mu(x)=1$.
It turns out
that this is the only source of non-uniqueness.

\begin{corollary}\label{Cor 2}
Assume that $\mu=1$ and that
 we are given the boundary $\p M$ and the response operator $\Lambda$.
Then using the
the processed time reversal iteration
we can find constructively the manifold $(M,g)$ upto an isometry
and on it the operator $A$ uniquely.

Moreover, if $M\subset \R^m$, $m\geq 2$, and $A$ is of the form
(\ref{eq: example A}), given the set $M$ and the response operator
$\Lambda$ we can determine
$c(x)$ and $q(x)$ uniquely in $M$.
 \end{corollary}

We note that the unique determination of $(M,g)$ and $A$
is well known, see e.g.\ \cite{BK,KK,KKL} and references therein.
The novelty of Corollary \ref{Cor 2} is that in the
construction based on processed time reversal iteration
the use of iterated measurement avoid many computationally
demanding steps used required in \cite{KKL}.
An example of such step is a Gram-Schmidt orthogonalisation of
basis of spaces $L^2(\Gamma\times [t_1,t_2])$, $\Gamma\subset \p M$
with respect to an inner
product determined by $\Lambda_{2T}$. This is a typical
step used in traditional Boundary Control method
that is  very sensitive for measurement errors.

\section{Proof of the convergence of the iteration}
\label{sec:localization proofs}

\subsection{Controllability results.}
The seminal result implying controllability is Tataru's
unique continuation result (see \cite{Ta1,Ta3})
\begin {theorem} [Tataru]
\label {th:3.4}
Let $u$ be a solution
of  wave equation
\ba
 u_{tt}(x,t)+Au(x,t)=0.
\ea
 Assume
that
\beq
u|_{\Gamma\times (0,2\tau)}=0,\hbox{ } \p_\nu u|_{\Gamma
\times
(0,2\tau)}=0
\label {3.20}
\eeq
where $\Gamma\subset \p M, \, \Gamma \neq \emptyset$ is open,
and $\tau>0$.
Then,
\bfo
u(\bx,\tau)=0,\  \p_tu(\bx,\tau)=0\quad \hbox{ for }\bx\in M(\Gamma,\tau).
\efo
\end {theorem}

This result yields the following Tataru's controllability result, see e.g.\
\cite{KKL} and references therein.

\begin {theorem}
\label {th:3.3}
Let $\Gamma \subset \p M$ be open and $\tau>0$. Then the linear subspace
\bfo
\{ u^f(\tau)\in L^2(M(\Gamma,\tau)):
       \ f \in C^{\infty}_0(\Gamma\times [0,\tau]) \}
\efo
is dense in $L^2(M(\Gamma,\tau))$.
\end {theorem}

\subsection{Inner products.}

Let us recall the  Blagovestchenskii identity:
\begin {lemma}
\label {l:4.10}
Let $f,h\in L^2 (\p M\times [0,2T])$.
Then
\beq
\int_M u^f(x,T) {u^h(x,T)}\,dV_\mu(x)=
\label {4.60}
\eeq
\bfo
\int_{[0,2T]^2}  \int_{\p M} J(t,s) \big[f(t)  {(\Lambda_{2T} h)(s)}-
(\Lambda_{2T} f)(t)  {h(s)}\big]\,dS_g(x)dtds,
\efo
where $J(t,s)=\frac 12\chi_{L}(s,t)$, see (\ref{Eq: L}).
\smallskip

\end {lemma}

\noindent
The proof in slightly different context
 is given e.g.\ in \cite{KKL}.

\noindent {\bf Proof.}  Let
$
w(t,s)=\int_M u^f(t) {u^h(s)}\,dV_\mu.
$
Integrating by parts, we see that
\begin{align}
\nonumber
(\p^2_t-\p^2_s)w(t,s)
\hspace{-5em}\\
&= \nonumber
-\int_M \big[Au^f(t) {u^h(s)}- u^f(t)
{Au^h(s)}\big]\,dV_\mu(x)
\\ \label{right hand side of wave eq}
&=-\int_{\p M} \big[\p_\nu u^f(t)u^h(s)- u^f(t)
{\p_\nu u^h(s)}\big]\,dS_g
\\ \nonumber
&=\int_{\p M} \big[(-\p_\nu u^f(t)+\eta u^f(t)) u^h(s)- u^f(t)
(-\p_\nu u^h(s)+\eta u^h(s))\big]\,dS_g
\\ \nonumber
&= \int_{\p M} \big[f(t) \HL h(s)-
\HL f(t)
{h(s)}\big]\,dS_g.
\end{align}
Moreover,
\bfo
\left. w \right| _{t=0}=\left. w\right| _{s=0}=0, \quad
\left. \p _t w\right| _{t=0}= \left. \p_s w\right| _{s=0}=0.
\efo
Thus we can consider (\ref{right hand side of wave eq}) as
one dimensional wave equation with known right hand side
and vanishing initial and boundary data. Solving this initial
boundary value problem we see the claim.

%
\hfill\proofbox

Consider now the operator $K$ appearing in formula (\ref {4.60b})
in more detail.
The Schwartz kernel of $\Lambda_{2T}$
is the Dirichlet boundary value of the  Green's function
$G(x,x',t-t')$ satisfying
\begin{align}\label{eq: Wave green}
( \p_t^2+A)G_{x',t'}(x,t) &= \delta_{x'}(x)\delta(t-t')\quad \hbox{ in
} M\times \R_+,\\
G_{x',t'}|_{t=0} &= 0,\quad \p_tG_{x',t'}|_{t=0}=0,  \nonumber \\
B_{\nu,\eta} G_{x',t'}|_{\p M\times \R_+} &= 0, \nonumber
\end{align}
where $G_{x',t'}(x,t)=G(x,x',t-t')$.
Then
\ba
G(x,x',t-t')=G(x',x,t-t')
\ea
and we see that
\ba
\Lambda_{2T}^*=R_{2T} \Lambda_{2T} R_{2T}
\ea
where $R_{2T}f(x,t)=f(x,2T-t)$ is the time reversal map.

Thus, as $\Lambda_{2T}^*=R_{2T}\Lambda_{2T} R_{2T}$, we can write result of Lemma \ref {l:4.10} as
\beq\label{kaava A}
\int_M u^f(x, T) {u^h(x, T)}\,dV_\mu(x)
=\int_{\p M\times [0,2T]} (Kf)(x,t)\,h(x,t)\,dS_g(x)dt
\eeq
where $K$ is defined in (\ref{K-definition}).

\subsection{Proof of Theorem \ref{Main F}.}

To produce focusing waves,
let us consider the problem
\beq\label{eq: minimize 0}
\min_{h\in L^2(B)} \|u^f(T)-u^h(T)\|_{L^2(M,dV_\mu)}^2
\eeq
where $B=\bigcup_{j=1}^J (\Gamma_j\times [T-T_j,T]),$
 $\Gamma_j\subset \p M$ are open and $0\leq T_j\leq T$.

As the solution of this minimisation problem does not always exist,
we can regularise the problem and study
\beq\label{eq: minimize}
\min_{h\in L^2(\p M\times [0,2T])}
F(h,\alpha)
\eeq
where $\alpha\in (0,1)$ and
\ba
F(h,\alpha)=\bra  K(Ph-f),Ph-f\cet_{L^2(\p M\times
[0,2T],dS_g)}+\alpha\|h\|_{L^2(\p M\times [0,2T],dS_g)}^2.
\ea
We recall that  $P=P_B$ is multiplication with the characteristic function
of $B$, that is,
$
(P_Bh)(x,t)=\chi_{B}(x,t)\,\cdotp h(x,t).
$

By (\ref{kaava A}),
\ba
F(h,\alpha)=\|u^f(T)-u^{Ph}(T)\|_{L^2(M,dV_\mu)}^2+\alpha\|h\|_{L^2(\p M\times [0,2T])}^2.
\ea

The minimisation problem  (\ref{eq: minimize})
is equivalent to (\ref{eq: minimize 0})
when $\alpha=0$.

\begin{lemma} \label{lem: 7}
For given $\alpha\in (0,1)$ the problem
(\ref{eq: minimize}) has a unique minimiser. Moreover,
the minimiser is the unique solution of the equation
\beq\label{main}
(PKP+\alpha)h=PK f.
\eeq
\end{lemma}

\noindent {\bf Proof.} The minimisation problem is strictly convex, and the
map $h\mapsto u^h(T)$ is continuous $L^2(\p M\times [0,2T])\to
H^{s_1}(M)$, $s_1<5/6$ by \cite{Lasiecka2}, see also \cite{Lasiecka}.
Now if $h_j\to h$ weakly
in $L^2(\p M\times [0,2T])$ as $j\to \infty$
then
 $Ph_j\to Ph$ weakly in  $L^2(\p M\times [0,2T])$.
As the embedding $H^{s_1}(M)\to  L^2(M)$
is compact and
the map $h\mapsto u^h(T)$ is continuous
 $L^2(\p M\times [0,2T])\to H^{s_1}(M)$,
then
 $u^{Ph_j}(T)\to u^{Ph}(T)$
in $L^2(M)$.
Therefore,
\begin{align*}
\lim_{j\to\infty}\bra & K(Ph_j-f),Ph_j-f \cet_{L^2(\p M\times
[0,2T])}
\\
&=
\lim_{j\to\infty} \|u^{Ph_j}(T)-u^f(T)\|_{L^2(M,dV_\mu)}^2=
 \|u^{Ph}(T)-u^f(T)\|_{L^2(M,dV_\mu)}^2.
\end{align*}

Let $h_j$ be a sequence such that
$\lim_{j\to \infty} F(h_j,\alpha)=\inf_h F(h ,\alpha)$.
As $F(h,\alpha)\geq \alpha \|h\|_{L^2}^2,$ the sequence  $(h_j)$ is bounded.
Since $\|h\|_{L^2} \leq \liminf_{j\to \infty} \|h_j\|_{L^2}$, this implies that $F(h, \alpha) \leq \lim \inf F(h_j, \alpha)$.

 Thus
by choosing a subsequence of $h_j$, we can assume that
$h_j$ converge to $h$ in the weak topology of $L^2(\p M\times [0,2T])$,
and the limit $h$ is a global minimiser of $F(\cdotp,\alpha)$.

By computing
the Fr\'{e}chet derivative of the functional  $F(h, \alpha)$
in any direction $\tilde h\in L^2(\p M\times [0,2T])$
at a minimiser $h$, we see that
\begin{align*}
0 &= D_h(\bra K(Ph-f),Ph-f\cet+\alpha \|h\|^2_{L^2})\tilde h\\
 &= \bra \tilde h, P^*K(Ph-f) \cet+
\bra \tilde h,P^*K^*(Ph-f)\cet+
2\alpha\bra \tilde h, h\cet
\end{align*}
for all $\tilde h$.
Since $K^*=K$ and
$P^*=P$, this implies that  $h$ satisfies the equation
(\ref{main}).
As $K$ is non-negative, $PKP+\alpha I\geq \alpha I$
so that equation (\ref{main}) has a unique solution providing
 the minimiser $h=h(\alpha)$.
\hfill\proofbox

\medskip
Note that  equation (\ref{main})
 also implies that
\beq
\label{21'}
h(\alpha) \in \hbox{Ran} (P).
\eeq

Next we want to solve equation (\ref{main}) using iteration.
To this end, let $\omega\in \R_+$ be a constant such that
\ba
\omega>2(1+\|P KP\|).
\ea
Then the equation (\ref{main}) can be written as
\beq\label{main 2}
(I-S)h=\frac 1 \omega \, PKf,\quad\hbox{where}\quad
S=I-\frac {\alpha+ P KP}\omega.
\eeq
Since $PKP$ is non-negative and
\ba
\alpha I\leq \alpha I+ P KP\leq \frac \omega 2 I,
\ea
we see that $\|S\| < 1$.
Thus we can solve $h$ using iteration:
Let
\ba
F:=\frac 1 \omega P Kf = \frac 1\omega P(J\Lambda_{2T} -R \Lambda_{2T} R J)f,
\ea
$h_0=0$, and consider the iteration
\ba
h_{n+1}=Sh_n+F,\quad n=0,1,\dots.
\ea
As $h_n=Ph_n$ by (\ref{21'}), we can
write
the iteration in the form (\ref{eq: iteration 1}), and
$\lim_{n\to \infty}h_n=h(\alpha)$ in $L^2(\p M\times [0,2T])$.

Applying this algorithm can find the minimisers
$h=h(\alpha)\in \hbox{Ran}\,(P)$ of problem (\ref{main}). They converge
by the following lemma:

\begin{lemma} \label{lem: 8}We have
\ba
\lim_{\alpha\to 0} u^{h(\alpha)}(x,T)=\chi_{N}(x) u^f(x,T)
\ea
in $L^2(M)$, with $N$ as in Theorem~\ref{Main F}.
\end{lemma}

\noindent {\bf Proof.}
Clearly,
\beq
\nonumber
& &F(h,\alpha)=F_0+\|\chi_{N}(u^{Ph}(T)-u^f(T))\|_{L^2(M)}^2
+\alpha\|h\|_{L^2(\p M\times [0,2T])}^2,\label{F-estimate}\\
& &\quad \quad F_0=\| (1-\chi_{N})u^f(T)\|_{L^2(M)}^2.
\eeq
Here the first term does not depend on $h$.
By Theorem \ref{th:3.3}, for any $\e>0$
there is a $h_\e\in L^2(B)$ such that
\beq
\label{22'}
\| u^{Ph_\e}(T)-\chi_{N}u^f(T)\|_{L^2(M)}^2<\frac \e 2.
\eeq
Thus we have
\ba
F(h_\e,\alpha)=F_0+\frac \e2+\alpha\|h(\e)\|_{L^2(\p M\times [0,2T])}^2.
\ea
This shows that if $\alpha<\alpha(\e)$,
where
\ba
\alpha(\e) :=\frac \e {2\|h(\e)\|_{L^2(\p M\times [0,2T])}^{2}},
\ea
then
\ba
F(h_\e,\alpha)\leq F_0+ \e.
\ea
Thus the minimiser $h(\alpha)$ of $F(h,\alpha)$ with $\alpha<\alpha(\e)$
satisfies
\ba
F(h(\alpha),\alpha)\leq F_0+\e.
\ea
As by (\ref{21'}), $h(\alpha)=P h(\alpha)$
so that $h(\alpha)\in L^2(\Gamma\times [0,T])$, we get from (\ref{22'}) that,
for $\alpha <\alpha(\e)$,
\beq\label{F-estimate 2}
\|\chi_{N}(u^{h(\alpha)}(T)-u^f(T))\|_{L^2(M)}^2\leq  \e/2.
\eeq
As $\supp(u^{h(\alpha)}(T))\subset N$, the claim follows.
\hfill\proofbox
\medskip

Theorem \ref{Main F} follows from Lemmata
\ref{lem: 7} and \ref{lem: 8}.\hfill\proofbox

\section{Proofs for the focusing of the waves}\label{sec:focus}

Next we prove  Corollary \ref{Cor 1}.
\vspace{1ex plus 1ex}

\noindent {\bf Proof of Corollary~\ref{Cor 1}.} Let $(z(x),s(x))$ be the boundary normal
coordinates of $x\in M$, that is, $s(x)=d(x,\p M)$
and $z(x)$ is the closest point of $\p M$ to $x$ when such
a point is unique. When a  closest boundary point is not unique,
the boundary normal coordinates are not defined.

We consider first the claim of the corollary in the case
when $\widehat T<\tau(\widehat z)$. Then the boundary normal
coordinates near $\widehat x = \gamma_{\widehat z, \nu}(\widehat T)$ are
well defined and the metric tensor in these coordinates
 has the form
\beq\label{bc kaava}
g=\left(\begin{array}{cc} 1& 0\\ 0& [g_{\a\beta}(z,s)]_{ \a,\beta=1}^{m-1}
\end{array}\right),
\eeq
where $[ g_{\a\beta}(z,s)]\in \R^{(m-1)\times (m-1)}$
 is a smooth positive definite
matrix valued function near $(\widehat z, \widehat T)=(z(\widehat x), s(\widehat x))$. Then there
is a $C_1=C_1(\widehat z, \widehat T)>0$ such that
\ba
C_1^{-1}I\leq [ g_{\a\beta}(z,s)]\leq
C_1I
\ea
 near $(\widehat z, \widehat T)$. Using this we see that
there is a $C_2>0$ such that
\ba
C_2^{-1}\leq
\frac {\vol_g \big(M(\{\widehat z\},\widehat T)\setminus M(\p M,T_0) \big)}
{(\widehat T-T_0)(\widehat T^2-T_0^2)^{(m-1)/2}}
\leq C_2.
\ea
As the solution $u^f(x,t)$ is smooth for $f\in C^\infty_0(\p M\times \R_+)$,
\ba
\lim_{T_0\to \widehat T}
\frac {\chi_{M(\{\widehat z\},\widehat T)\cup M(\p M,T_0)}-\chi_{M(\p M,T_0)}}{\vol_g\big(M(\{\widehat z\},\widehat T)\setminus M(\p M,T_0)\big)}
\, u^f(T)=u^f(\widehat x,T)\delta_{\widehat x}
\ea
in ${\mathcal D}'(M)$, the claim follows in the case $\widehat T<\tau(\widehat z)$.

When  $\widehat T>\tau(z)$, the claim follows from the fact
that for sufficiently small $\widehat T-T_0>0$
the set $M(\{\widehat z\},\widehat T)\setminus M(\p M,T_0)$ is empty.
\hfill\proofbox

\medskip

\noindent
{\bf Remark}\, The above proof  of
 Corollary \ref{Cor 1} also shows that in  Corollary \ref{Cor 1}
\beq\label{c kaava}
C_0(x_0)=\lim_{T_0\to \widehat T} \frac
{(\widehat T-T_0)^{(m+1)/2}}
{\vol_g\big(M(\{\widehat z\},\widehat T)\setminus M(\p M,T_0)\big)}.
\eeq

\section{Boundary distance functions and reconstruction of the metric}\label{sec:metric}

Here we present a method for finding the boundary distance
functions using processed time reversal iteration.

Let $z,y\in \p M$, $0\leq T_1\leq \tau(z)$ and
$T>\dist_{\p M}(y,z)+T_1$.
 Denote
$x=\gamma_{z,\nu}(T_1)$. Next we give an algorithm
that can be used to determine $\dist(x,z)$.

To this end, let
$\Gamma_j\subset \p M$
and  $\Sigma_j\subset \p M$
be neighbourhoods of $z$ and $y$, respectively, such
that
 $\Gamma_j\to \{z\}$ and
$\Sigma_j\to \{y\}$ when $j\to \infty$.

Let $\e>0$, $\tau\in [0,T]$, and
\begin{align*}
 N_j^1 &= M(\Gamma_j,T_1)   & B_j^1    &= \Gamma_j\times [T-T_1,T]\\
 N_j^2 &= M(\Sigma_j,\tau)  & B_j^2    &= \Sigma_j\times [T-\tau,T]\\
 N_{\e}^3 &= M(\p M,T_1-\e) & B_{\e}^3 &= \p M\times [T-(T_1-\e),T].
\end{align*}
\begin{lemma}\label{lem: test part 1}
The distance  $\dist(x,z)$ is the infimum of all $\tau\in [0,T]$ that satisfy
the condition
\beq\label{cond: 1}
& &\hbox{the set $I(j,\e)=(N_j^1\cap N_j^2)\setminus N_{\e}^3$
contains} \\ \nonumber
& & \hbox{a non-empty open set  for all }
j\in \Z_+,\ \e>0.
\eeq
\end{lemma}

\begin{figure}
  \centering
  \psfrag{1}{$x$}
  \psfrag{2}{$z$}
  \psfrag{3}{$y$}
  \psfrag{4}{$\tau$}
  \psfrag{5}{$T_1$}
  \psfrag{6}{$\Gamma_j$}
  \psfrag{7}[c]{$\Sigma_j$}
  \psfrag{A}{$T_1-\e$}
  \includegraphics[width=6cm]{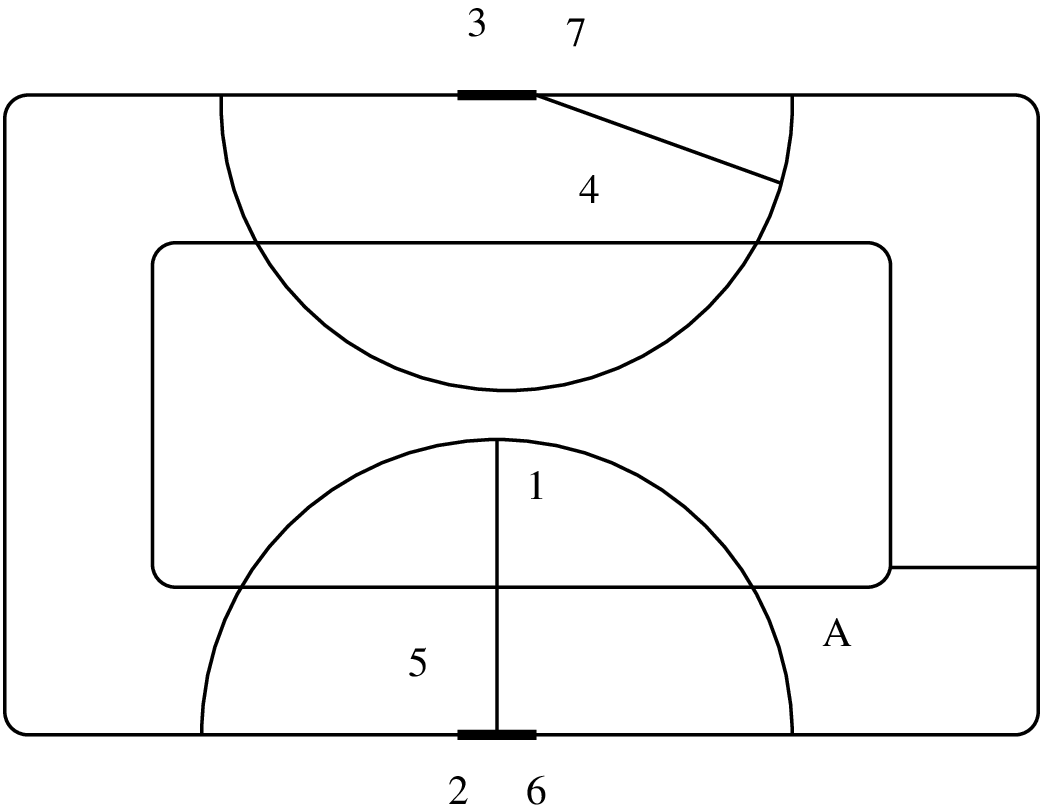} \hspace{.6cm}\includegraphics[width=6cm]{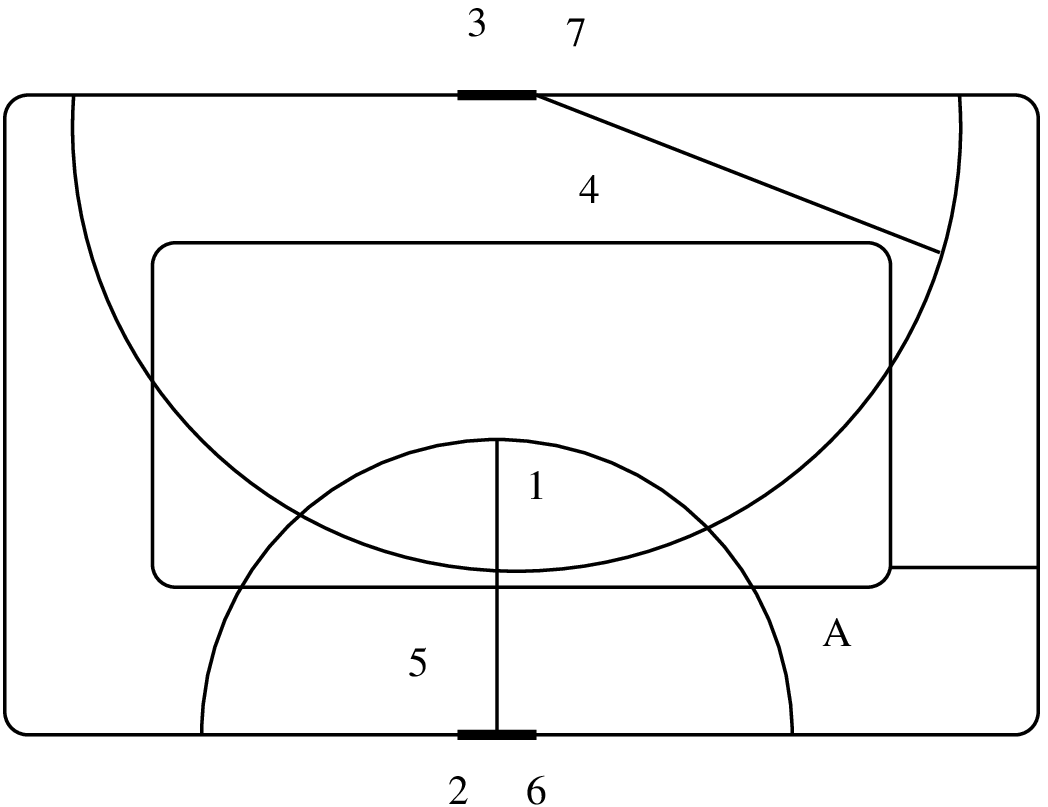}
  \caption{Sketch of the situation in Lemma~\ref{lem: test part 1}
    when $I(j, \epsilon)$ is empty (left) and non-empty (right).}
  \label{fig:lem test part 1 fig}
\end{figure}

See Figure~\ref{fig:lem test part 1 fig} for a sketch of the cases
of empty and non-empty $I(j,\e)$.

\vspace{1ex plus 1ex}
\noindent {\bf Proof.}
First consider what happens  if $\dist(x,y)<\tau$. Since
$T_1\leq \tau(z)$,  we see that then
$B(x,r)\cap (N_j^1\setminus N_{\e}^3)$ contains a non-empty open set
for all $r>0$, where $B(x,r)\subset M$ is a ball of $(M,g)$
with center $x$ and radius $r$.
When $r<\tau-\dist(x,y) $, we see that $B(x,r)\subset N_j^2$. Thus
$I(j,\e)$ contains an open set and (\ref{cond: 1}) is satisfied.

On other hand,
if $\dist(x,y)>\tau$, let $r=\dist(x,y)-\tau$. When $j\to \infty$
 and $\e\to 0$, we see using metric in the boundary
normal coordinates (\ref{bc kaava}) that $N_j^1\setminus N_{\e}^3\to
\{x\}$ in the Hausdorff metric. Thus when $j$ is
large enough and $\e$ is small enough,
$N_j^1\setminus N_{\e}^3\subset B(x,r/2)$. Then
$ B(x,r/2)\cap N_j^2=\emptyset$, and
(\ref{cond: 1}) is not satisfied.

Summarizing, condition (\ref{cond: 1}) is satisfied if $\dist(x,y)<\tau$
and not satisfied if $\dist(x,y)>\tau$.
As $\dist(x,y)<T$, this yields the claim.
\hfill\proofbox
\medskip

Next we show that by using the processed time reversal iteration,
we can test if condition (\ref{cond: 1}) is valid. To this end,
denote
\begin{align*}
\tilde N_1(j, \e) &=N^1_{j}\cup N_{\e}^3 &\tilde B_1(j, \e) &=B^1_{j}\cup B_{\e}^3\\
\tilde N_2(j, \e) &=N^2_{j}\cup N_{\e}^3 &\tilde B_2(j, \e) &=B^2_{j}\cup B_{\e}^3\\
\tilde N_3(j, \e) &=N^1_{j}\cup N^2_{j}\cup N_{\e}^3
     &\tilde B_3(j, \e) &=B^1_{j}\cup B^2_{j}\cup B_{\e}^3\\
\tilde N_4(j, \e) &=N_{\e}^3 &\tilde B_4(j,\e) &=B_{\e}^3.
\end{align*}
Let   $f\in C^\infty_0(\p M\times \R_+)$.
Using the processed time reversed iteration on time interval
$[0,2T]$ with projectors $P_B$
corresponding to $B=\tilde B_k(j, \e)$, $k=1,2,3,4$ and starting point $f$,
 we obtain functions
$
h_n(\alpha;\e,j,k)\in L^2(\p M\times [0,2T]).
$
Using them, define
\beq\label{p-functions}
 p_n(\alpha,j,\e)=h_n(\alpha;\e,j,1)+h_n(\alpha;\e,j,2)-h_n(\alpha;\e,j,3)
-h_n(\alpha;\e,j,4).
\eeq

\begin{lemma}\label{lem: test part 2}
The condition (\ref{cond: 1}) is satisfied if and only if
there exists an $f\in C^\infty_0(\p M\times \R_+)$ such that
for any $j\in \Z_+$ and $\e>0$ the functions $p_n(\alpha,j,\e)$ defined
in formula (\ref{p-functions}) satisfy
\beq\label{cond 2}
\lim_{\alpha\to 0}\lim_{n\to \infty}\bra K_{2T}
p _n(\alpha,j,\e),p_n(\alpha,j,\e)\cet\not =0.
\eeq
\end{lemma}

\noindent
{\bf Proof.} The functions
$h_n(\alpha;\e,j,k)$ defined by processed time reversal iteration satisfy
\ba
\chi_{\tilde N_k(j,\e)}u^f(T)=\lim_{\alpha\to 0} \lim_{n\to \infty}
u^{h_n(\alpha;\e,j,k)}(T),\quad k=1,2,3,4.
\ea
A simple computation gives us
\beq\label{chi formula}
\chi_{I(j,\e)}(x)=\chi_{\tilde N_1(j,\e)}(x)+
\chi_{\tilde N_2(j,\e)}(x)-\chi_{\tilde N_3(j,\e)}(x)-\chi_{\tilde N_4(j,\e)}(x)
\eeq
for all $x\in M$.
Therefore, using (\ref{p-functions}) we see that  in $L^2(M)$
\ba
\chi_{I(j,\e)}u^f(T)
=\lim_{\alpha\to 0}\lim_{n\to \infty} u^{p_n(\alpha,j,\e)}(T).
\ea
By Theorem \ref{th:3.3} we see that the functions
$u^f(T)$ with $f\in C^\infty_0(\p M\times \R_+)$ are
smooth functions that form a dense
set in $L^2(M(\p M,T))$.
Thus condition (\ref{cond: 1}) is satisfied if and only if
there exists an
$f\in C^\infty_0(\p M\times \R_+)$ such that for any $j\in \Z_+$ and $\e>0$
\ba
 & &\bra\chi_{I(j,\e)}u^f(T),u^f(T)
\cet_{L^2(M)}\\
& &=\lim_{\alpha\to 0}\lim_{n\to \infty}\bra
K_{2T}
p _n(\alpha,j,\e),p_n(\alpha,j,\e)\cet_{L^2(\p M\times [0,2T])}\not =0.
\ea
This proves the claim.
\hfill\proofbox

We note that if condition (\ref{cond 2}) is satisfied
for some $f\in C^\infty_0(\p M\times \R_+)$, it is satisfied
for all such $f$ in an open and dense set.
\medskip

Lemmata \ref{lem: test part 1} and \ref{lem: test part 2}
give an algorithm for the determination of $d(y,x)$ by using
processed time reversal iteration
\begin{align}
\label{algorithm}
d(y,x)&=
\inf\{\tau\in  [0,T]:\ \hbox{there is an $
f\in C^\infty_0(\p M\times \R_+)$}
\\
\nonumber &\quad\quad\quad\quad \hbox {
 such that (\ref{cond 2}) holds for all $j\in \Z_+$ and $\e>0$}\}.
\end{align}

Summarizing, we have proven:

\begin{proposition}\label{prop: boundary distannce} Assume we are given $\p M$ and the response
operator $\Lambda$. Let $z,y\in \p M$, $T_1\leq \tau(z)$.
Then using the algorithm (\ref{algorithm})
we can compute $\dist(x,y)$ for
$x=\gamma_{z,\nu}(T_1)$.
\end{proposition}

Let us consider consequences of the result above. To this end,
we define the set of the boundary
distance functions. For each $x\in M$,
the corresponding boundary
distance function, $r_x \in C(\p M)$ is  given by
\[
 r_x: \p M\to \R_+,\quad r_x(z)=\dist(x,z), \quad
z \in \partial M.
\]
In fact, $r_x\in \hbox{Lip}\,(\p M)$ with the Lipschitz constant equal to one.
The boundary
distance functions define {\it the
 boundary distance map} ${\mathcal R}:M\to C(\p M)$,
${\mathcal R}(x)=r_x$, which is continuous and injective
(see \cite {KKL}).  Denote by
\[
 {\mathcal R}(M)=\{r_x\in C(\partial M): \, x\in M\},
\]
the image
of ${\mathcal R}$. It is known (see \cite{KKL,KKLima})
  that, given the set
 ${\mathcal R}(M)
\subset C(\partial M)$
we can endow it, in a constructive way,
 with
a  differentiable structure and a metric tensor $\tilde g$,
so that $({\mathcal R}(M),\tilde g)$ becomes a manifold that is
 isometric
to $(M,g)$,
\[
 ({\mathcal R}(M),\tilde g)\cong (M,g).
\]
See \cite{KKLima}. A stable construction procedure of $(M,g)$ as
a metric space from the set ${\mathcal R}(M),\tilde g)$
and H\"older type stability estimates are give in
\cite{KatsudaKurylevLassas}.
\medskip

\noindent{\bf Example}\, By the triangle inequality,
\beq\label{triangle}
\|r_x-r_y\|_\infty\leq \dist (x,y).
\eeq
We consider in this example the case when $(M,g)$ is a compact
manifold such that all points $x,y\in M$ can be joined with a
unique  shortest geodesic. This implies that, for any
$x,y\in M$,  the shortest geodesic $\gamma([0,s])$ from $x$ to $y$, parametrised
along arc length, can be continued to a maximal
geodesic $\gamma([0,L])$ that hits  the boundary at
a point $z=\gamma(L)\in \p M$. Then
$$|r_x(z)-r_y(z)|=\dist (x,y)$$ implying equality in
(\ref{triangle}). Thus in the case when all geodesics between
arbitrary points $x,y\in M$ are  unique, the manifold $(M,g)$
is isometric to the manifold ${\mathcal R}(M)$ with the distance function
inherited from $C(\p M)$. In the general case the construction
of the metric is more elaborate.

\medskip

By Theorem \ref{prop: boundary distannce} we can compute for
all $x=\gamma_{z,\nu}(T)$ with $T\leq \tau(z)$ the corresponding
boundary distance
function $r_x$. Since all points $x\in M$ can be represented
in this form (see e.g.\ \cite{chavel})
 we  can find the set ${\mathcal R}(M)$
that can be endowed with a manifold structure
isometric to the original manifold $(M,g)$.
We have thus proven the following result:

\begin{corollary}\label{cor; manifold} Assume we are given $\p M$ and the response
operator $\Lambda$. Then using
the  processed time reversal iteration
we can find the manifold $(M,g)$ upto an isometry.
\end{corollary}

Corollary \ref{cor; manifold} can also be formulated by
saying that we can find the metric tensor $g$ in local coordinates.
For example, for any point $x_0\in M$ there are $z_j\in \p M$, $j=1,2,\dots,m$,
 such that
$x\mapsto (X^1(x),\dots,X^m(x))$ with $X^j(x):=\dist(x,z_j)$ define local coordinates
near $x_0$. In these coordinates the distance functions $x\mapsto\dist(x,z)$, $z\in \p M$
determine the metric tensor. For details of this construction, see \cite{KKL}.

Next we prove Corollary \ref{Cor 2}.
\vspace{1ex plus 1ex}

\noindent {\bf Proof of Corollary~\ref{Cor 2}.}
As the boundary data $\p M$ and $\Lambda$ determine
$(M,g)$ upto isometry, we can apply formula (\ref{c kaava}) to
find the function $C_0(x)$, $x\in M$
in Corollary \ref{Cor 1} in any local coordinates of $(M,g)$. Thus
in local coordinates we can find the values
of waves $u^f(x,t)$ for all $x\in M$, $t>0$. By Tataru's controllability
theorem, the waves $u^f(x,T)$ with fixed $T>0$
form a dense set in $L^2(M(\p M,T))$.
In this dense set we can find in local coordinates the distributions
\ba
A u^f(x,t)=-\p_t^2u^f(x,t)= -u^{f_{tt}}(x,t),
\ea
implying that we can find values $Aw$ for all $w\in L^2(M(\p M,T))$.
As $T$ is arbitrary, we find  $Aw$ for all $w\in L^2(M)$.

The case when $M\subset \R^m$, $m\geq 2$ and $A$ is of the form
(\ref{eq: example A}), we use the fact
that a manifold with conformally Euclidian metric is
embedded into $\R^m$ in a unique way when $\p M$ is fixed.
Thus we find $c(x)$ when a metric $g$ is given, and as
above we also find
an operator $A_\kappa$ that is equivalent to $A$ through a gauge
transformation. As there is a unique operator of the form
(\ref{eq: example A}) in the orbit of the
the gauge transformations $A\to A_\kappa$ the claim follows
(see \cite{KKL}).
\hfill\proofbox

\section{Iteration when measurements have errors}

Let $(\Omega,\Sigma,\prob)$ be a complete probability space.

Assume the measurements have a random noise, that is,
measurements give us with an input $f$ the function
$\Lambda_{2T} f+\e$, where $\e$ is random Gaussian noise
that has values in $L^2(\p M\times (0,T))$. Assume that
$ \expec \e=0$ and denote the covariance operator of $\e$ by $C_\e$.
Note that $C_\e$ is a compact operator on $L^2(\p M\times (0,T))$
(see \cite{bogachev})
and thus the standard white noise on $\p M\times (0,T)$ does not
satisfy our assumptions.

Assume that the noise is independent of previous
measurements each time when we do a new measurement.
When the noise is added to  the processed time reversal iteration,
we come to the iteration of the form
\ba
\tilde h_{n+1}=S\tilde h_n+p+N_n,
\ea
where $N_n=PJ\e_n^1-PR\e_n^2$ where $\e_n^1$ and $\e_n^2$ are random variables
having the same distribution
as $\e$. Thus the $N_n$ are independent identically distributed
Gaussian random variables satisfying
${\Bbb E}N_n=0$
and having covariance operator $C_N=PJC_\e J^*P^*+PRC_\e R^*P^*$.
Let us consider the averaged results of iterations
\beq\label{random limit}
 \tilde h_K^{ave}= \frac 1K\sum_{n=1}^K \tilde h_n.
\eeq
Then
\begin{align*}
 \frac 1K\sum_{n=1}^K  \tilde h_{n}
&=
\frac 1K\sum_{n=1}^K  h_{n}+
\frac 1K\sum_{n=1}^K \sum_{m=1}^n (n-m+1)S^{n-m}N_n\\
&=
\left(\frac 1K\sum_{n=1}^K h_{n}\right)+
\left(\frac 1K\sum_{n=1}^K (I-S)^{-2}N_n\right)
\\
& \qquad + \left(
\frac 1K\sum_{n=1}^K
(-(I-S)^{-2}S^{n+2}- (n+2)
(I-S)^{-1}S^{n+1})N_n\right)\\
&=H^1_K+H^2_K+H^3_K,
\end{align*}
where $h_n$ are the results of the processed
time reversal iteration without noise.

Above, the deterministic term $H^1_K$ converges to the same limit
$\lim_{n\to \infty} \tilde h_n=h(\alpha)$
as the processed time reversal iteration without noise.
Now consider $H^2_K$ and $H^3_K$ as random variables in
$L^2(\p M\times (0,T))$. They can also be viewed as
random fields on $\p M\times (0,T)$, see \cite{Rozanov}.
By the law of large numbers in infinite dimensional spaces
\cite{HJ-P}, 
 we see that
\beq\label{eq: large numbers}
\lim_{K\to \infty}\|H^2_K\|_{L^2(\p M\times (0,T))}=0.
\eeq
As $\|S\|_{L(L^2(\p M))}\leq \frac 12$, the last term $H^3_K$
satisfies also an estimate  analogous to (\ref{eq: large numbers}).
Thus  the averaged processed time reversal iteration with noise
converges to the same limit as processed time reversal iteration
without noise, that is,
\ba
\lim_{K\to \infty} \tilde h_K^{ave}=
h(\alpha)\quad \hbox{in }
L^1(\Omega;L^2(\p M\times (0,T))).
\ea

\medskip

\noindent
{\bf Acknowledgements:} The research was partially supported by the
TEKES project ``MASIT03 -- Inverse Problems and Reliability of Models''.

\bigskip

{\obeylines
Kenrick Bingham and Matti Lassas
Institute of Mathematics
P.O.Box 1100
02015 Helsinki University of Technology
Finland 
}
\bigskip

{\obeylines
Yaroslav Kurylev
Department of Mathematical Sciences  
Loughborough University
Loughborough
Leicestershire
LE11 3TU  
UK
}
\bigskip

{\obeylines
Samuli Siltanen
Institute of Mathematics
Tampere University of Technology
P.O.Box 553 
33101 Tampere
Finland
}

\end{document}